\documentclass[12pt]{article}

\usepackage{psfig,amssymb}

\makeatletter
\def\eqalign#1{\null\vcenter{\def\\{\cr}\openup\jot\m@th
  \ialign{\strut$\displaystyle{##}$\hfil&$\displaystyle{{}##}$\hfil
      \crcr#1\crcr}}\,}
\makeatother

\newcommand{\be}{\begin{equation}} 
\newcommand{\ee}{\end{equation}}
\newcommand{\beq}{\begin{eqnarray}}
\newcommand{\eeq}{\end{eqnarray}}
\newcommand{\bt}{\begin{theorem}}
\newcommand{\et}{\end{theorem}}
\newcommand{\bl}{\begin{lemma}}
\newcommand{\el}{\end{lemma}}
\newcommand{\bc}{\begin{corollary}}
\newcommand{\ec}{\end{corollary}}
\newcommand{\ba}{\begin{array}}
\newcommand{\ea}{\end{array}}

\newcommand{\la}{\label}
\newcommand{\ci}{\cite}

\newtheorem{theorem}{Theorem}
\newtheorem{lemma}[theorem]{LEMMA}
\newtheorem{corollary}[theorem]{COROLLARY}

\newcommand{\ti}{\tilde}
\newcommand{\de}{\delta}
\newcommand{\De}{\Delta}
\newcommand{\al}{\alpha}
\newcommand{\ga}{\gamma}

\newcommand{\si}{\sigma}
\newcommand{\Si}{\Sigma}
\newcommand{\om}{\omega}

\newcommand{\ze}{\zeta}
\renewcommand{\th}{\theta}

\newcommand{\ep}{\varepsilon }

\newcommand{\bi}{\bibitem}

\newfont{\msbm}{msbm10 scaled\magstep1}
\newfont{\msbms}{msbm7 scaled\magstep1} 

\newcommand{\bbc}{\mbox{$\mbox{\msbm C}$}}

\begin{document}
\def\wt{\widetilde}
\bigskip\bigskip\bigskip
\begin{center}
{\Large\bf 
Gap probability in the spectrum of random matrices
and asymptotics of polynomials orthogonal on an arc of the unit circle. 
}\\
\bigskip\bigskip\bigskip
I. V. Krasovsky\\
\bigskip Technische Universit\"at Berlin\   
Institut f\"ur Mathematik MA 7-2\\
Strasse des 17. Juni 136, D-10623, Berlin, Germany\\
E-mail: ivk@math.tu-berlin.de\\
\bigskip\bigskip
\end{center}
\bigskip\bigskip\bigskip
\noindent{\bf Abstract.} We obtain uniform asymptotics for polynomials
orthogonal on a fixed and varying arc of the unit circle with a 
positive analytic weight function. We also
complete the proof of the large $s$ asymptotic expansion for the
Fredholm determinant with the kernel $\sin z/(\pi z)$ on the interval 
$[0,s]$, verifying a conjecture of Dyson for the
constant term in the expansion.
In the Gaussian Unitary Ensemble of random matrices, this
determinant describes the probability for an interval of length $s$
in the bulk scaling limit to be free from the eigenvalues.
\bigskip\bigskip\bigskip

\section{Introduction}
One problem in the random matrix theory is estimation of the probability
for a given interval to be free from the eigenvalues. In the Gaussian
Unitary Ensemble this probability for any interval of length $2s$
in the bulk scaling limit is equal to the following
Fredholm determinant:
\begin{equation}
\De(s)=\det[I-K],
\label{De}
\end{equation}
where $K$ is the integral operator on $L^2(0,2s)$ given by
\[
(Kg)(x)=\int_0^{2s}{\sin (x-y)\over \pi(x-y)}g(y)dy.
\]

The probability of a gap in the spectrum of random matrices from
orthogonal and symplectic ensembles 
is also expressed in terms of $\De(s)$ (see \cite{Mehta,TW,Dbook}). 

It was shown by Jimbo,  Miwa, M\^ori, and Sato \cite{JMMS} 
that $(d/ds)\ln \De(s)$ satisfies a 
modified  Painlev\'e V equation (for simpler proofs of this see 
\cite{TW,DIZ}).  

An interesting question is 
calculation of the asymptotics of $\De(s)$ for large $s$ 
(the small $s$ series are easy to obtain). 
The first two terms in the expansion
of $\ln\De(s)$ were found by des Cloizeaux and Mehta \cite{CM} who used
a connection with the spheroidal functions. The full asymptotic
expansion was obtained by Dyson \cite{Dyson} with the help of the inverse
scattering techniques for Schr\"odinger operators.  
These calculations were partly conjectural. A rigorous derivation of the
main term was given by Widom \cite{WidomAs} using continuous analogues of
orthogonal polynomials. Finally, Deift, Its, and Zhou \cite{DIZ} 
(see also that work for a
more extensive introduction) calculated, as a particular case of a
more general result, the full asymptotics of the derivative
$(d/ds)\ln\De(s)$ using techniques of matrix Riemann-Hilbert problems.
This settled the question up to the constant term
in the expansion of $\ln\De(s)$.
The first 3 terms in the Dyson expansion are as follows:
\be
\ln \De(s)=-{s^2\over 2}-{1\over 4}\ln s+c_0+O\left({1\over s}\right),
\qquad s\to\infty,\label{exp}
\ee
where the constant term $c_0=(1/12)\ln 2 +3\zeta'(-1)$, and
$\zeta'(x)$ is the derivative of Riemann's zeta function.
Thus, justification of $c_0$ here remained the only problem and it is
solved in the present paper. Actually, we obtain the first 
3 terms:\footnote{As this paper was being prepared for publication,
an announcement by T. Ehrhardt claiming the same result as Theorem 1
(by a different method) was posted on the internet.

A third solution to the problem by a Riemann-Hilbert approach (related to
the present one) is in preparation by P. Deift, A. Its, and X. Zhou.
}
\bt
The large $s$ asymptotics of $\ln \De(s)$ are given by (\ref{exp}).
\et

The proof is based on a formula by Deift \ci{D} which connects the
determinants of two Toeplitz matrices, a formula by Widom \ci{Warc} for
asymptotics of Toeplitz determinants on a circular arc, and on
asymptotics for orthogonal polynomials on a circular arc which are
computed here.

Let $f_\al(\th)$ be a weight function on an arc 
$\al\le\th\le 2\pi-\al$, $z=e^{i\th}$, 
$0<\al<\pi$ of the unit circle $|z|=1$, 
and $\phi_n(z,\al)=\chi_n z^n+\cdots$, $n=0,1,\dots$ 
the corresponding system of orthonormal polynomials:
\be
{1\over 2\pi}\int_\al^{2\pi-\al}
\phi_k(e^{i\th},\al)\overline{\phi_m(e^{i\th},\al)}f_\al(\th)d\th=
\de_{km},\qquad k,m=0,1,\dots.\la{phi1}
\ee

Such polynomials in the case of the circle ($\al$=0) were first 
studied by Szeg\H o (see \ci{Sz}) who, in
particular, found several important asymptotics for them as $n\to\infty$.  
Afterwards, asymptotic analysis of such polynomials was carried
out by many authors. Specifically for the case of an arc (whose study
was initiated by Akhieser \ci{Akh}), see
\ci{GN,GAkh,BL,RL} and references therein. 
However, the full asymptotic expansion
at all points $z\in\bbc$ for a wide class of weights
became a feasible task only after recent development
of Riemann-Hilbert problem methods. 
It was observed by Fokas, Its, and Kitaev \ci{fik} that 
orthogonal polynomials satisfy certain matrix Riemann-Hilbert problems.
An efficient method for their asymptotic solution (steepest descent
techniques) was developed by Deift and Zhou \ci{dz,Dbook}
and applied for analysis of polynomials orthogonal on the real axis
in \ci{DK1,DK2} (see also \ci{BI} for a different Riemann-Hilbert
approach) and on the unit circle in \ci{bdj}. 
The case of polynomials orthogonal on $[-1,1]$
(especially relevant for the present work) 
was considered by Kuijlaars, McLauphlin, Van Assche, and Vanlessen
\ci{kuimcl,kuivan} who found full asymptotics at all points in the
case of a positive analytic weight on $[-1,1]$ with 
power-type singularities at the end-points. 

In the present paper, we shall give a procedure to obtain full
asymptotics for all $z$ for polynomials $\phi_n(z,\al)$ and their leading
coefficients as $n\to\infty$ in the case of a positive analytic
weight $f_\al(\th)$. 
The argument will be similar to that of  \ci{kuimcl}.
We consider $2s/n\le\al<\pi$, $n>s$, $s\to\infty$, which includes 
both the cases of a fixed arc ($\al$ is
independent of $n$) and a varying arc.
The asymptotics for $\phi_n(z,\al)$ we obtain are in the inverse
powers of $n\sin(\al/2)$. 
The remainder after $k$ terms is uniform in $\al$.
The general solution is given by equations
(\ref{i2}--\ref{aschi}). The first 2 asymptotic terms for any $z$ can
be easily written using (\ref{R1},\ref{A1B1}). An example is given in
(\ref{p1},\ref{c1}). 

Our solution can be generalized to the following cases:
(1) the weight $f_\al(\th)$ has power-type singularities at the
end-points of the arc (this can be done following \ci{kuimcl});
(2) the weight $f_\al(\th)$ is not analytic but only smooth enough and
positive (one can approximate it then by its Fourier series).

If the weight is symmetric $f_\al(\th)=f_\al(2\pi-\th)$, there exist
relatively simple formulas of Szeg\H o type \ci{Sz,Zh,fredh}
connecting polynomials on an arc with those on an interval. In this
case and for a fixed arc one could try to obtain our results
from those of \ci{kuimcl, kuivan}. 
The present argument, however, is more direct.

For the proof of Theorem 1, we shall only need asymptotics of
polynomials with the weight $f_\al(\th)=1$ and only at the
point $e^{i\al}$. What we need is summarized in the following theorem
proved in Section 2 (after the argument in the general case is given):

\bt
Let $0<\al<\pi$, $\ga=\cos(\al/2)$,
\be\eqalign{
r^1_+={e^{-i\al/2}\over 3\cdot2^4 i}(1+e^{-i\al}-2e^{i\al}),\qquad
r^1_-={e^{-i\al/2} \over 2^4 i}(1+3e^{i\al}),\\
r^2_+={1\over 3\cdot2^9}(16-9e^{i\al}+43e^{-i\al}-2e^{-2i\al}),\qquad
r^2_-={1\over 2^9}(-6+7e^{i\al}-17e^{-i\al}),\\
{r^1_-}'={e^{i\al/2} \over 4}\cos^2(\al/2),\qquad
\tau={1+2\cos\al\over 6i},\qquad
\rho=n\sin(\al/2),\qquad \ep>0.}
\ee
Let $f_\al(\th)=1$.
Then the polynomial $\phi_n(z,\al)$ admits an asymptotic expansion 
for large $\rho$ in the inverse powers of $\rho$. We have for $z=e^{i\al}$
\be
\phi_n(e^{i\al},\al)=\chi_n
\ga^n e^{i\al(n/2-1/4)}\sqrt{\pi i \rho}
\left[1+{r^1_- \over \rho}+
{r^2_- \over \rho^2}+ {r^3_- \over \rho^3}+
O\left({1\over\rho^4}\right)\right],\la{phin}
\ee
where $\chi_n$ is the leading coefficient of 
$\phi_n(z,\al)=\chi_n z^n+\cdots$ for which we have
\be
\chi_{n-1}^2=\ga^{-2n+1}\left[1+{1\over 4n}+{5\over
    2^5n^2}+O\left({1\over n^3}\right)\right],\la{chin}
\ee
and $r^3_-$ is a bounded function of $\al$.
The derivative of the polynomial
$\phi_n'(z,\al)=(d/dz)\phi_n(z,\al)$
at $e^{i\al}$ can be written as
\be\eqalign{
\phi_n'(e^{i\al},\al)={n\over 2}\phi_n(e^{i\al},\al)e^{-i\al}+
\chi_n\ga^n e^{i\al(n/2-5/4)}{\sqrt{\pi i \rho}\over 2\sin\al}
\left[i\rho^2+e^{i\al/2}\rho+\tau+\right. \\ \left.
{1\over\rho}(r^1_-(i \rho^2+\tau)+r^1_+e^{i\al/2}\rho+{r^1_-}')+
{1\over\rho^2}(r^2_- i\rho^2+r^2_+e^{i\al/2}\rho)+
{i r^3_-\over\rho} +O\left({1\over\rho^2}\right)\right].}\la{phinp}
\ee
There exists $s_0>0$ such that
all the remainder terms are valid and uniform 
in $\al$, $s$, and $n$
for $\al\in[2s/n,\pi-\ep]$, $s>s_0$, and $n>s$.
\et
 
\noindent
{\bf Remark}
The uniformity of the remainders here is crucial for 
the proof of Theorem 1.

After constructing asymptotics for polynomials in Section 2, we give a
proof of Theorem 1 in Section 3. Note that the present method could
also be used to obtain the full asymptotic expansion of $\De(s)$.

\section{Asymptotics of polynomials on an arc}
In the present section we construct asymptotics for polynomials
$\phi_n(z)=\chi_n z^n+\cdots$ orthonormal with 
a weight $f(z)=f_\al(\th)$ on an arc
$\al\le\th\le2\pi-\al$, $z=e^{i\th}$, for $2s/n\le\al\le\pi-\ep$, $\ep>0$, 
$n>s$, $s\to\infty$. This includes both the cases of a fixed arc
$0<\al<\pi$ and the varying arc $\al=2s/n$. The function $f(z)$ is
assumed positive and analytic on the arc for a fixed arc case, and on
the whole circle in the general case. In the general case, we obtain
$\al$-uniform asymptotics in the inverse powers of $n\sin(\al/2)$. 

Consider the following $2\times 2$ matrix 
\begin{equation} \label{RHM}
    Y(z) =
    \pmatrix{
\chi_n^{-1}\phi_n(z) & 
\chi_n^{-1}\int_{\Si}{\phi_n(\xi)\over \xi-z}
{f(\xi)d\xi \over 2\pi i \xi^n} \cr
\chi_{n-1}\phi^*_{n-1}(z) & 
\chi_{n-1}\int_{\Si}{\phi^*_{n-1}(\xi)\over \xi-z}
{f(\xi)d\xi \over 2\pi i \xi^n}
    },
\end{equation}
where $\Si$ is the arc $\al \le\th\le 2\pi-\al$ of the unit circle
traversed in the direction from $2\pi-\al$ to $\al$, and
$\phi^*_n(z)=z^n\overline{\phi_n(1/{\bar z})}$.
As is easy to verify (see \ci{bdj,kuimcl}),
$Y(z)$ is the unique solution of the following Riemann-Hilbert problem:
\begin{enumerate}
    \item[(a)]
        $Y(z)$ is  analytic for $z\in\bbc \setminus\Si$.
    \item[(b)]
For $\th \in (\al,2\pi-\al)$, 
$Y$ has continuous boundary values
$Y_{+}(x)$ as $z$ approaches $x=e^{i\th}$ from
the outside of the circle, and $Y_{-}(x)$, from the inside. 
They are related by the jump condition
\begin{equation}\label{RHPYb}
            Y_+(x) = Y_-(x)
            \pmatrix{
                1 & x^{-n}f(x) \cr
                0 & 1},
            \qquad\mbox{$x=e^{i\th},\qquad \th\in (\al,2\pi-\al)$.}
        \end{equation}
    \item[(c)]
        $Y(z)$ has the following asymptotic behavior at infinity:
        \begin{equation} \label{RHPYc}
            Y(z) = \left(I+ O \left( \frac{1}{z} \right)\right)
            \pmatrix{
                z^{n} & 0 \cr
                0 & z^{-n}}, \qquad \mbox{as $z\to\infty$.}
        \end{equation}
    \item[(d)]
       Near the end-points of the arc $e^{\pm i\al}$ 
        \begin{equation}\label{RHPYd}
            Y(z) =  O\pmatrix{ 
                    1 & \log|z-e^{\pm i\al}| \cr
                    1 & \log|z-e^{\pm i\al}|},
        \end{equation}
as $z \to e^{\pm i\al}$, $z \in \bbc \setminus \Si$.
\end{enumerate}

The function
\be
\psi(z)={1\over2\ga}\left(z+1+\sqrt{(z-e^{i\al})(z-e^{-i\al})}\right),
\qquad \ga=\cos(\al/2),\la{psi}
\ee
which conformally maps the outside of the arc $\Si$ into the outside
of the unit circle, will have an 
important role in what follows (cf. (6.18) of \ci{DIZ}). 
Here we take the branch of the square
root which is positive for positive arguments.
Note that the boundary values of $\psi(z)$, $\psi_+(x)$ as $z$
approaches $x\in\Si$ from the outside of the circle, and  $\psi_-(x)$, from
the inside, are related as:
\be
\psi_+(x)\psi_-(x)= x.\la{bvf}
\ee
Let $\mu(z)$ be defined by the equation (\ref{psi}) but with the minus
sign in front of the square root. Then $\mu(z)$ is the mapping of the
outside of the arc into the {\it inside} of the unit circle. 
Hence $|\mu(z)|<1$, whereas $|\psi(z)|>1$ for 
$z\in\bbc\setminus\Si$.
Therefore we have that first, for $|z|<1$  $|z/\psi(z)^2|<1$, and
second, for $|z|>1$
\[
\left|{z\over\psi(z)^2}\right|=
\left|{z\mu(z)^2\over\psi(z)^2\mu(z)^2}\right|=
\left|{z\mu(z)^2\over z^2}\right|=\left|{\mu(z)^2\over z}\right|<1.
\] 
Thus,
\be\la{ineq}
\left|{z\over\psi(z)^2}\right|<1
\qquad\mathrm{for}\qquad |z|\neq 1.
\ee
This inequality will be useful later on.

We now replace the original Riemann-Hilbert problem with an equivalent one
which is normalized to unity at infinity and has oscillating elements
of the jump matrix.
Namely, set
\be
T(z)=\ga^{-n\si_3}Y(z)\psi(z)^{-n\si_3},\qquad \si_3=\pmatrix{1&0\cr 0&-1}.
\ee
Then, as is easy to verify, $T(z)$ satisfies the same problem
as $Y(z)$ (\ref{RHM}) but with the changed conditions (b) and (c):
\begin{enumerate}
\item[(b)]
\begin{equation}
            T_+(x) = T_-(x)
            \pmatrix{
                x^n\psi_+(x)^{-2n} & f(x) \cr
                0 & x^n\psi_-(x)^{-2n}},
            \qquad\mbox{$x\in\Si$,}
        \end{equation}
    \item[(c)]
        \begin{equation} 
            T(z) = I+ O \left( \frac{1}{z} \right), 
     \qquad \mbox{as $z\to\infty$.}
\end{equation}
\end{enumerate}

Following the idea of the steepest descent method of Deift and Zhou, 
we now replace the Riemann-Hilbert problem for $T(z)$ with 
an equivalent one on a system of 3 contours where some of the 
jump matrix elements are exponentially small.
Divide the complex plane into 3 regions as shown on Figure 1
(the contours $\Si_{1,3}$ lie sufficiently close to $\Si_2\equiv\Si$
for $f(z)$ to remain nonzero and analytic in regions 1 and 2)
\begin{figure}
\centerline{\psfig{file=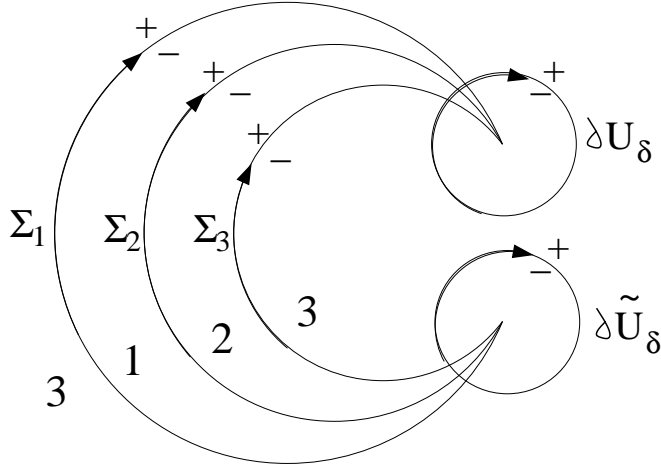,width=4.0in,angle=-90}}
\vspace{0.1cm}
\caption{
The contour for Riemann-Hilbert problems.}
\label{fig1}
\end{figure}
and define the matrix-valued function $S(z)$ by the formulas:

\noindent
1) in region 1
\be
S(z)=T(z)\pmatrix{1&0\cr -{z^n\over f(z)\psi(z)^{2n}}& 1},
\ee 
2) in region 2
\be
S(z)=T(z)\pmatrix{1&0\cr {z^n\over f(z)\psi(z)^{2n}}& 1},
\ee
3) in region 3
\be
S(z)=T(z).
\ee

The condition (c) in the problem for $S(z)$ is the same as
for $T(z)$, the conditions (a), (b) and (d) are different. Namely,
\begin{enumerate}
\item[(a,b)]
$S(z)$ is analytic in $\bbc\setminus(\Si_1\cup\Si_2\cup\Si_3)$ with the
following jump conditions on the contours:
\begin{equation}
            S_+(x) = S_-(x)
            \pmatrix{1&0\cr {x^n\over f(x)\psi(x)^{2n}}& 1},
\qquad\mbox{$x \in \Si_1\cup\Si_3$,}
\end{equation}
\begin{equation}
S_+(x) = S_-(x)
            \pmatrix{0&f(x)\cr
              -f(x)^{-1}&0},
\qquad\mbox{$x \in \Si_2\equiv\Si$,}
\end{equation}
\item[(c)]
as $z\to\infty$
\be
 S(z) = I+ O \left( \frac{1}{z} \right), 
\ee
\item[(d)]
near the end-points of the arc
\begin{equation}
            S(z) =  O\pmatrix{ 
                   \log|z-e^{\pm i\al}| & \log|z-e^{\pm i\al}| \cr
                   \log|z-e^{\pm i\al}| & \log|z-e^{\pm i\al}|},\la{scS}
        \end{equation}
as $z \to e^{\pm i\al}$, $z \in \bbc \setminus (\Si_1\cup\Si_2\cup\Si_3)$.
\end{enumerate}
Recalling (\ref{ineq}),
we see that, for $n$ large, for $x$ outside some neighborhoods of the
endpoints of the arc, the jump matrix on $\Si_1\cup\Si_3$ is uniformly
exponentially close to the identity. We therefore approximate the
function $S(z)$ with parametrices inside the mentioned neighborhoods
and in the outside region where we neglect the jumps on
$\Si_1\cup\Si_3$. The Riemann-Hilbert problems for the parametrices
can be solved and the solution closely resemble that in case of
polynomials on an interval \ci{kuimcl}.

The parametrix for the outside region is defined as the solution of the
following Riemann-Hilbert problem:
\begin{enumerate}
    \item[(a)]
        $N(z)$ is  analytic for $z\in\bbc \setminus\Si\equiv\Si_2$,
    \item[(b)]
with the jump condition on $\Si$
\be
N_+(x) = N_-(x)
            \pmatrix{0&f(x)\cr
              -f(x)^{-1}&0},
\qquad\mbox{$x \in \Si$},
\ee
\item[(c)]
and the following behavior at infinity
\be
N(z) = I+ O \left( \frac{1}{z} \right), 
     \qquad \mbox{as $z\to\infty$.}
\ee
\end{enumerate}
The solution $N(z)$ is found in the same way as in \ci{kuimcl}.
Consider the Szeg\H o function:
\be
{\cal D}(z)=\exp\left[{\sqrt{(z-e^{i\al})(z-e^{-i\al})}\over 2\pi i}\int_\Si
{\ln f(\xi)\over \sqrt{(\xi-e^{i\al})(\xi-e^{-i\al})}}
{d\xi \over \xi-z}\right].\la{Sf}
\ee
This function is analytic outside the arc, and its boundary values
satisfy
\be
{\cal D}_+(x){\cal D}_-(x)=f(x),\qquad x\in\Si.
\ee
Denote
\be
{\cal D}_\infty=\lim_{z\to\infty} {\cal D}(z)=
\exp\left[-{1\over 2\pi i}\int_\Si{\ln f(\xi) d\xi \over 
\sqrt{(\xi-e^{i\al})(\xi-e^{-i\al})}}\right].\la{Sfinf}
\ee

Then the solution of the above Riemann-Hilbert problem is as follows:
\be
N(z)={1\over2}({\cal D}_\infty)^{\si_3}
\pmatrix{a+a^{-1}&-i(a-a^{-1})\cr i(a-a^{-1})& a+a^{-1}}
{\cal D}(z)^{-\si_3},
\qquad
a(z)=\left({z-e^{i\al}\over z-e^{-i\al}}\right)^{1/4},\la{N}
\ee
where the value of the root satisfies the condition $a(z)\to 1$ as
$z\to\infty$. Note that $\det N(z)=1$, which allows, in particular, to
write a simple expression for the inverse $N(z)^{-1}$.

Now consider a $\de$-neighborhood $U_\de$ of the point $e^{i\al}$,
small enough so that $f(z)$ is analytic and nonzero there.
The jump
matrices on $\Si_1$, $\Si_3$ are not close to the identity in this
region, so we need to construct a separate local parametrix.
We look for a matrix-valued function $P(z)$ which is analytic in
$U_\de$, satisfies the same jump relations on 
$(\Si_1\cup\Si_2\cup\Si_3)\cap U_\de$ as $S(z)$, has the same
behavior at $z=e^{i\al}$ as $S(z)$, and matches $N(z)$ at the
boundary:
\be
P(z)N(z)^{-1}=I+O\left({1\over n\sin(\al/2)}\right),\qquad z\in 
\partial U_\de\setminus (\Si_1\cup\Si_2\cup\Si_3),\la{match}
\ee
where $\rho=n\sin(\al/2)\to\infty$.
We seek $P(z)$ in the form
\be
P(z)=E(z)\hat P(z)\left({\psi(z)\over\sqrt{z}}\right)^{-n\si_3}
f(z)^{-\si_3/2},\la{P1}
\ee
where $E(z)$ is invertible and analytic in a neighborhood of $U_\de$.
The function $E(z)$ does not affect jump relations and will 
be chosen later so that $P(z)$ satisfies (\ref{match}).
Using the boundary-value property (\ref{bvf}), we obtain 
(cf. \ci{kuimcl}):
\be\eqalign{
\hat P(x)_+=\hat P(x)_-
\pmatrix{1&0\cr 1&1}, \qquad \mbox{$x\in (\Si_1\cup\Si_3)\cap U_\de,$}\\
\hat P(x)_+=\hat P(x)_-
\pmatrix{0&1\cr -1&0}, \qquad \mbox{$x\in \Si_2\cap U_\de.$}}\la{jtP}
\ee
Consider the function $\om(z)$ defined by the equation:
\be
e^{\sqrt{\om(z)}}={\psi(z)\over\sqrt{z}}\la{g}
\ee
Using (\ref{bvf}), we have for the boundary values of $\om$ on the arc:
\be
\sqrt{\om(x)}_+=\ln{\psi_+(x)\over\sqrt{x}}=-\sqrt{\om(x)}_-,
\ee
therefore $\om(z)$ is analytic in $U_\de$. 
For $z$ near $e^{i\al}$, we obtain uniformly for all $\al$
\be\eqalign{
{\psi(z)\over\sqrt{z}}=
1+\left(\sqrt{i\sin(\al/2)\over\cos(\al/2)}e^{-i\al/2}
\sqrt{1+{u\over 2i\sin\al}}\sqrt{u}+
(e^{i\al/2}\cos(\al/2)^{-1}-1)e^{-i\al}{u\over 2}\right)(1+O(u)),
\\
u=z-e^{i\al}.}\la{psi_as}
\ee
If also $|u|<2\sin\al$, we get for the function $\om(z)$
the following (nonuniform) expansion at $z=e^{i\al}$:
\be
\om(z)={i\sin(\al/2)\over \cos
  (\al/2)}e^{-i\al}u\left(
1-{1-2e^{-i\al}-2e^{-2i\al}\over 6i\sin\al}u+O(u^2)\right).\la{g_as1}
\ee
Denote
\be
\hat P(z)=Q(\zeta),\qquad \zeta=n^2 \om(z).
\ee
Now we reached a crucial moment.
The circle $\partial U_\de$ is transformed in
the $\zeta$ variable into a curve $\partial\hat U_\de$ 
whose minimal distance from zero is
$n^2 \min_{0\le t\le 2\pi}|\om(e^{i\al}+\de e^{it})|$.
In order to construct a solution, we need this distance to be large.
This is so for large $n$ if $\de$ and $\al$ 
are independent of $n$ (see (\ref{g_as1})).
In the general case of $2s/n\le\al\le\pi-\ep$, there exists some 
small $\al_0$ (depending on $f(z)$, $\ep$) such that we can assume 
\be
\de=\cases{\sin(\al_0/2),& for $\al_0\le\al\le\pi-\ep$\cr
\sin(\al/2),& for $2s/n\le\al\le\al_0$}.\la{del}
\ee
Putting $u=\de e^{it}$ in (\ref{psi_as}) and choosing $\al_0$
sufficiently small, we obtain after simple analysis:
\be\la{g_as2}
n^2\min_{0\le t\le 2\pi} |\om(e^{i\al}+\de e^{it})|\ge C(n\sin(\al/2))^2
\ee
for some constant $C$ which is larger than zero. We see that for large
$s$ this distance remains uniformly large (not less than of order
$s^2$) for any $\al\in[2s/n,\pi-\ep]$, where $n>s$.

In $\hat U_\de$ the image $\hat\Si_{1,2,3}$ of the cuts can be 
considered as 3 direct
lines emanating from zero. (The image of $\Si$ is a line, and the exact
form of $\Si_1$ and $\Si_3$ can be chosen at will.) 
The analytic matrix-valued function $Q(\zeta)$ which satisfies the jump
conditions (\ref{jtP}) on $\hat\Si_{1,2,3}$ and singularity conditions
(\ref{scS}) at $\zeta=0$ was constructed in \ci{kuimcl}.
Namely, we have in $\hat U_\de$ (the regions 1, 2, and 3 correspond to the
$\zeta$-variable images of the regions in Figure 1) the following
expressions in terms of modified Bessel and Hankel functions (see,
e.g., \ci{Abr}):

1) region 1
\begin{equation}\label{Q1}
    Q(\zeta) ={1\over 2}
    \pmatrix{
        H_{0}^{(1)}(e^{-i\pi/2}\zeta^{1/2}) &
        H_{0}^{(2)}(e^{-i\pi/2}\zeta^{1/2}) \cr
        \pi \zeta^{1/2} \left(H_{0}^{(1)}\right)'(e^{-i\pi/2}\zeta^{1/2}) &
        \pi \zeta^{1/2} \left(H_{0}^{(2)}\right)'(e^{-i\pi/2}\zeta^{1/2})
    },
\end{equation}

2) region 2 
\begin{equation}\label{Q2}
    Q(\zeta) ={1\over 2}
    \pmatrix{
        H_{0}^{(2)}(e^{i\pi/2}\zeta^{1/2}) &
        -H_{0}^{(1)}(e^{i\pi/2}\zeta^{1/2}) \cr
        -\pi \zeta^{1/2} \left(H_{0}^{(2)}\right)'(e^{i\pi/2}\zeta^{1/2}) &
        \pi \zeta^{1/2} \left(H_{0}^{(1)}\right)'(e^{i\pi/2}\zeta^{1/2})
    },
\end{equation}

3) region 3
\begin{equation}\label{Q3}
    Q(\zeta) =
    \pmatrix{
     I_{0} (\zeta^{1/2}) & \frac{i}{\pi} K_{0}(
        \zeta^{1/2})\cr
     \pi i \zeta^{1/2} I_{0}'(\zeta^{1/2}) & 
     -\zeta^{1/2} K_{0}'(\zeta^{1/2})},
\end{equation}
where $-\pi<\arg(\zeta)<\pi$.

We now have to choose $E(z)$ so that the matching condition
(\ref{match}) is satisfied. For that we can use the first term in 
the asymptotic expansion of Bessel and Hankel functions
for large $\zeta$. 
The expansion for $Q(\ze)$ is the same in all the 3 regions. We can
write down an arbitrary number of terms. Below we shall make use only
of the first three. We have:
\be\eqalign{
Q(\ze)={1\over\sqrt{2}}(\pi\sqrt\ze)^{-\si_3/2}\pmatrix{1&i\cr i&1}
\left[I+{1\over 8\sqrt\ze}\pmatrix{-1&-2i\cr -2i&1}-\right.\\ \left.
{3\over 2^7\ze}\pmatrix{1&-4i\cr 4i&1}+
O(\ze^{-3/2})\right]e^{\sqrt\ze\si_3}}\la{Qas}
\ee
uniformly on the boundary of $\hat U_\de$.
We now define $E(z)$ as follows:
\be
 E(z)={1\over\sqrt{2}}N(z)f(z)^{\si_3/2}
\pmatrix{1&-i\cr-i&1}(\pi n\sqrt{\om(z)})^{\si_3/2}.\la{E}
\ee
We verify exactly as in \ci{kuimcl} that it is an analytic function 
in $U_\de$.
Using (\ref{P1}), (\ref{g}), and (\ref{g_as2})
(and estimating $N(z)$ as below)
we see that the matching 
condition (\ref{match}) is now satisfied.
Thus we have
\be
P(z)=E(z)Q(n^2 \om(z))\left({\psi(z)\over\sqrt{z}}\right)^{-n\si_3}
f(z)^{-\si_3/2}.\la{P}
\ee
Solution in the neighborhood $\ti U_\de$ of $e^{i(2\pi-\al)}$ 
is similar (but
note the reversed direction of the contours). We have there:
\be\eqalign{
\ti P(z)=\ti E(z)\si_3 Q(n^2 \om(z))\si_3
\left({\psi(z)\over\sqrt{z}}\right)^{-n\si_3}
f(z)^{-\si_3/2},\\
\ti E(z)={1\over\sqrt{2}}N(z)f(z)^{\si_3/2}
\pmatrix{1&i\cr i&1}(\pi n\sqrt{\om(z)})^{\si_3/2}.}\la{P2}
\ee
We are now ready for the last transformation of the Riemann-Hilbert
problem. Let
\be\eqalign{
R(z)=S(z)N^{-1}(z),\qquad
z\in\bbc\setminus(\overline{U_\de\cup\ti U_\de}\cup\Si_{1,2,3}),\\
R(z)=S(z)P^{-1}(z),\qquad
z\in U_\de\setminus\Si_{1,2,3},\\
R(z)=S(z)\ti P^{-1}(z),\qquad
z\in \ti U_\de\setminus\Si_{1,2,3}.}\la{Rb}
\ee
It is easy to see that this function has jumps only on 
$\partial U_\de$, $\partial\ti U_\de$, and parts of $\Si_1$, and
$\Si_3$ lying outside of the neighborhoods $U_\de$, $\ti U_\de$ (we
denote these parts $\Si_{1,3}^\mathrm{out}$).
Namely,
\be\eqalign{
R_+(x)=R_-(x)N_-(x)\pmatrix{1&0\cr {x^n\over f(x)\psi(x)^{2n}}&1}
N_-(x)^{-1},\qquad x\in\Si_{1,3}^{\mathrm out},\\
R_+(x)=R_-(x)P(x)N(x)^{-1},\qquad x\in \partial U_\de,\\
R_+(x)=R_-(x)\ti P(x)N(x)^{-1},\qquad x\in \partial\ti U_\de.}\la{Rj}
\ee
As is easy to verify, the
matrix elements of $N(x)$ and $N(z)^{-1}$ remain bounded 
for $|x-e^{\pm  i\al}|\ge\de\ge\sin(s/n)$ and all $\al$.
For example,
\be\la{a4}
|a^4|=(1+8\ga\sin t+16\ga^2)^{-1/2}\le(1-4\ga)^{-1}
\qquad\mathrm{for}\qquad z=e^{i\al}+\sin(s/n)e^{it}\in\partial U_{\sin(s/n)}.
\ee
(For evaluation of ${\cal D}(z)$ near $e^{\pm i\al}$ see a similar
calculation in Lemma 6.4 of \ci{kuimcl}.)
Therefore the jump matrix on $\Si_{1,3}^\mathrm{out}$
can be uniformly (both in $z$ and $\al$)
estimated by (\ref{g},\ref{g_as2}) as
$I+O(\exp{(-C_1 n\sin(\al/2))})$, where $C_1$
is a positive constant. Note that in the case of a fixed arc, it is
sufficient to use (\ref{ineq}).

The jump matrices on $\partial U_\de\cup\partial\ti U_\de$
admit an asymptotic expansion in powers of $1/\sqrt{\ze}$ (which turns
into expansion in powers of
$1/n$ for a fixed arc or $1/s$ for a varying arc case). 
First, (\ref{Qas}), (\ref{E}), and (\ref{P}) yield
\be\eqalign{
P(z)N(z)^{-1}=
I+N(z)f(z)^{\si_3/2}\left[
{1\over 8\sqrt\ze}\pmatrix{-1&-2i\cr -2i& 1}-
{3\over 2^7\ze}\pmatrix{1&-4i\cr 4i& 1}+\right.\\ \left.
O(\ze^{-3/2})\right]
f(z)^{-\si_3/2}N(z)^{-1}=
I+\De_1+\De_2+O(\ze^{-3/2}),\qquad z\in \partial U_\de,}\la{lastjump}
\ee
where $\De_1$ and $\De_2$ denote the terms with $\sqrt\ze$ and $\ze$,
respectively (the remainder term will be justified below).
We write them down explicitly for the case of $f(z)=1$
needed later on:
\be
\De_1={1\over 2^4\sqrt\ze}\pmatrix{
-(3a^2-a^{-2})& -i(3a^2+a^{-2})\cr
-i(3a^2+a^{-2})& 3a^2-a^{-2}},
\qquad
\De_2=-{3\over 2^7\ze}\pmatrix{1& -4i\cr 4i& 1},\la{de12}
\ee
where $f=1$. The functions $\De_1(z)$
and $\De_2(z)$ for $z\in\partial\ti U_\de$ are given by the same
expressions with $a$ exchanged with $1/a$ and $i$ replaced by $-i$.
 
Using the expansion for Bessel functions, we can write a general term
$\De_j$ in (\ref{lastjump}) which is of order $1/\ze^{j/2}$. 
Indeed,
apart from the prefactor with $\ze^{j/2}$, the matrix elements
of $\De_j$ are obviously $O(1)$ as $\ze\to\infty$ for a fixed arc.
It is also true in the general case, because of 
the remark after equation (\ref{Rj}).
Now it is clear that (\ref{lastjump}) is an asymptotic expansion in $\ze$.

Since for $z\in\partial U_\de$ by (\ref{g_as1},\ref{g_as2})
\be
\sqrt\ze=n\sqrt{\om(z)}=O(n\sin(\al/2))=
\cases{O(n),& $\al$ fixed\cr
O(s),& $\al=2s/n$,} 
\ee
the component $\De_j$ on this boundary is of order $1/n^j$ (fixed arc)
and $1/s^j$ (varying arc), and 
the remainder term in (\ref{lastjump}) is
uniform for all $\al$ between a fixed positive value and $2s/n$,
all $n>s$, and all $s$ larger than some $s_0$. (As we shall see this
uniformity persists for the remainder in the asymptotics of our polynomials.)
We show as in \ci{kuimcl}
that $\De_j(z)$ is an analytic function in $U_\de$ with a 
pole at $e^{i\al}$ of order less than or equal to $[(j+1)/2]$.
The same reasoning also holds for the neighborhood $\ti U_\de$.
We shall denote the components of the jump matrix there by the same
symbols $\De_j(z)$.

If $R(z)$ is known, we can trace the sequence $Y\mapsto T\mapsto
S\mapsto R$ backwards, and obtain an expression for the polynomials. 
 
We look for $R(z)$ asymptotically in the 
form $R(z)\sim R_0(z)+R_1(z)+R_2(z)+\cdots$,
where $R_j(z)$ is of the same order as $\De_j$ (in our case, of order
$(n\sin(\al/2))^{-j}$).
For more discussion and
justification of this expansion see \ci{DK1,kuimcl}. More precisely,
it can be shown as in Theorems 7.8--7.10 of \ci{DK1} that
for any $k\ge1$
\be
R(z)=I+\sum_{j=1}^{k-1}R_j(z)+O((n\sin(\al/2))^{-k})\la{R}
\ee
uniformly for all $z$ if the arc is fixed, and for $z$ outside a
neighborhood of $z=1$ if the arc is varying ($\al=2s/n$).
The proof of this expansion in the general case $2s/n\le\al<\pi$
for $z$ close to
$1$ requires a special argument since the contours then are close to $z$
and can not be trivially deformed. Such a proof will be given later on
for the case of $f=1$, $z\in U_\de$, we need below.
The argument in the general case is similar.
Moreover, it follows directly from the proofs that the remainder 
term in (\ref{R}) has the same uniformity property
as that in (\ref{lastjump}).

Substituting this asymptotic expansion
into (\ref{Rb}) and collecting the terms of the same order, we obtain:
\be\eqalign{
R_{0+}(x)+R_{1+}(x)+\cdots\sim
(R_{0-}(x)+R_{1-}(x)+\cdots)(I+\De_1(x)+\cdots),
\qquad x\in \partial U_\de\cup\partial\ti U_\de.\\
R_{0+}(x)=R_{0-}(x)\quad\Rightarrow\quad R_0(z)=I,\\
R_{1+}(x)-R_{1-}(x)=\De_1(x),\\
R_{2+}(x)-R_{2-}(x)=R_{1-}(x)\De_1(x)+\De_2(x),\\
R_{k+}(x)-R_{k-}(x)=\sum_{j=1}^k R_{k-j,-}(x)\De_j(x),\qquad k=1,2,\dots}
\ee
The main term in the asymptotics of polynomials 
is given therefore by
the parametrices at the appropriate points $z$.
The expressions for $R_k(z)$ follow from the Sokhotsky-Plemelj
formulas:
\be
R_1(z)={1\over 2\pi i}\int_{\partial U_\de\cup\partial\ti U_\de}
{\De_1(x)dx\over x-z},\qquad
R_2(z)={1\over 2\pi i}\int_{\partial U_\de\cup\partial\ti U_\de}
{R_{1-}(x)\De_1(x)+\De_2(x)\over x-z}dx,\quad\dots\la{plem}
\ee
Note that the contours are traversed in the negative direction (see
Figure 1).

Following \ci{kuimcl}, we can also obtain the expressions for $R_j(z)$ 
in a different way.
As mentioned,  $\De_1(z)$ is analytic in $U_\de\cup\ti
U_\de$ with simple poles at the end-points of the arc. Thus,
\be
\De_1(z)=\frac{A^{(1)}}{z-e^{i\al}}+O(1),\quad\mbox{as } z\to e^{i\al}, 
\qquad
    \De_1(z)=\frac{B^{(1)}}{z-e^{-i\al}}+O(1),\quad\mbox{as } z\to
    e^{-i\al},
\ee 
where the constant matrices $A^{(1)}$ and $B^{(1)}$ are obtained by expanding
various functions in (\ref{lastjump}) at $z=e^{i\al}$ and 
$z=e^{-i\al}$. It is easy to verify directly that 
the Riemann-Hilbert problem for $R_1(z)$ has the solution: 
\be
R_1(z)=
\cases{\frac{A^{(1)}}{z-e^{i\al}}+\frac{B^{(1)}}{z-e^{-i\al}},&
for $z\in\bbc\setminus(\overline{U_\de\cup\ti U_\de})$\cr
\frac{A^{(1)}}{z-e^{i\al}}+\frac{B^{(1)}}{z-e^{-i\al}}-\De_1(z),&
for $z\in U_\de\cup\ti U_\de$.}\la{R1}
\ee

Expanding the functions in (\ref{lastjump}), we obtain:
\be
A^{(1)}={\cos(\al/2)\over 8n}\pmatrix{1&-i{\cal D}_\infty^2\cr
-i{\cal D}_\infty^{-2}& -1}e^{i\al/2},\qquad
B^{(1)}=\overline {A^{(1)}},\la{A1B1}
\ee
where $\overline {M}$ means complex conjugation applied
to every matrix element of $M$.
The general term $R_k(z)$ is obtained similarly provided we have 
the expressions for $R_j(z)$, $j=1,2,\dots,k-1$.  
Since $R_j(z)$ are analytic in $U_\de\cup\ti U_\de$
and $\De_j(z)$ have poles at $e^{\pm i\al}$ of order 
at most $[(j+1)/2]$, we see that
\be
\sum_{j=1}^k R_{k-j,-}(z)\De_j(z)=
\frac{A^{(k)}_p}{(z-e^{i\al})^p}+
\frac{A^{(k)}_{p-1}}{(z-e^{i\al})^{p-1}}+\cdots+
\frac{A^{(k)}_1}{z-e^{i\al}}+O(1),\quad \mbox{as } z\to e^{i\al},\la{i1}
\ee
where $p=[(k+1)/2]$, and similar expressions hold with matrices $A$
replaced with some matrices $B$ in a neighborhood of $e^{-i\al}$.
Then
\be
R_k(z)=
\cases{
\sum_{j=1}^p\left(\frac{A^{(k)}_j}{(z-e^{i\al})^j}+
\frac{B^{(k)}_j}{(z-e^{-i\al})^j}\right),&
for $z\in\bbc\setminus(\overline{U_\de\cup\ti U_\de})$\cr
\sum_{j=1}^p\left(\frac{A^{(k)}_j}{(z-e^{i\al})^j}+
\frac{B^{(k)}_j}{(z-e^{-i\al})^j}\right)-
\sum_{j=1}^k R_{k-j,-}(z)\De_j(z),&
for $z\in U_\de\cup\ti U_\de$}.\la{i2}
\ee

Recalling the definitions of $R$, $S$, and $T$, we finally get
\be
Y(z)=\ga^{n\si_3}R(z)M(z)\psi(z)^{n\si_3},\qquad
R(z)\sim I+R_1(z)+R_2(z)+\cdots,
\ee
where for $z$ restricted to region 3 $M(z)=N(z)$, $P(z)$, or $\ti P(z)$ if 
$z\in\bbc\setminus(\overline{U_\de\cup\ti U_\de})$, $U_\de$, or $\ti U_\de$,
respectively (the expressions in regions 1 and 2 can also be
readily written).
Therefore, by (\ref{RHM}),
\be
\chi_n^{-1}\phi_n(z)=Y_{11}(z)=\ga^n\psi^n(z) 
[R_{11}(z)M_{11}(z)+R_{12}(z)M_{21}(z)].\la{aspol}
\ee
Furthermore,
\be
\chi_{n-1}^2=Y_{21}(0)=\ga^{-2n}
[R_{21}(0)N_{11}(0)+R_{22}(0)N_{21}(0)].\la{aschi}
\ee
Since we know the expressions for $M(z)$ and can obtain $R(z)$ with an
arbitrary precision, the last 2 equations give an implicit solution
for asymptotics of polynomials $\phi_n(z)$ and their leading
coefficients. These are the asymptotic series in the
inverse powers of $n\sin(\al/2)$. The error
after $k$ terms is $O(n\sin(\al/2))^{-k-1}$
and remains uniform for all $\al$, $s$, and $n$, provided
$\al\in[2s/n,\pi-\ep]$, $n>s$, $s>s_0$. 

As an example, we give below the first two terms in the asymptotics of
$\phi_n(z)$ valid for $z$ outside a fixed arbitrary small
$\epsilon$-neighborhood of the arc $\Si$:
\be
\chi_n^{-1}\phi_n(z)\sim \ga^n\psi^n(z){{\cal D}_\infty\over {\cal
  D}(z)}\left[{a(z)+a(z)^{-1} \over 2}+{\ga\over 8n}
\left({a(z)e^{i\al/2}\over z-e^{i\al}}+
{a(z)^{-1}e^{-i\al/2}\over z-e^{-i\al}}\right)\right],\la{p1}
\ee
\be
\chi_{n-1}^2\sim {\ga^{-2n+1}\over {\cal D}(0){\cal D}_\infty}
\left[1+{1\over 4n}\right],\la{c1}
\ee
where all the quantities are defined as above (see (\ref{psi},
\ref{Sf}, \ref{Sfinf}, \ref{N})).

We now give an explicit solution for the first 3 terms in the 
asymptotics of $\phi_n(z)$
and its first derivative at an end-point of the arc for the weight $f(z)=1$. 

Recall that we still have to prove the remainder term in the expansion (\ref{R})
for $z\in U_\de$, $2s/n\le\al<\al_0$.
Since $f(z)=1$ we see that ${\cal D}(z)=1$.
Each matrix element of $\De_k(z)$ can be written as 
$(\al_k a^2+\beta_k a^{-2}+\ga_k)\ze^{-k/2}$, where $\al_k$, $\beta_k$,
and $\ga_k$ are independent of $\al$. In a neighborhood of $e^{i\al}$,
we have the series:
\be\eqalign{
\ze^{-k/2}=(n^2 u\sin\al)^{-k/2}\sum_{j=0}^\infty b_j(k,\al)
\left({u\over\sin\al}\right)^j,\qquad
a(z)^2=\sqrt{u\over\sin\al}
\sum_{j=0}^\infty c_j \left({u\over\sin\al}\right)^j,\\
|u|<\sin(\al/2),\qquad u=z-e^{i\al},
\qquad 0<|c_j|\le1,\quad j\ge0,}
\ee
where $b_j(k,\al)$ are bounded functions of $\al$.
Since $\De_k(z)$ is single-valued in $U_\de$, its matrix
elements do not contain terms with $\sqrt{u}$. Hence, 
$\al_k,\beta_k=0$ if $k$ is even, and $\ga_k=0$ if $k$ is odd.
We have in the same neighborhood:
\be
\De_k(z)={1\over (n\sin\al)^k}\sum_{j=-[(k+1)/2]}^\infty C_j(k,\al)
{u^j\over \sin^{j}\al}.
\ee
Now using (\ref{i1}) and the second formula in (\ref{i2}), 
it is easy to show by induction that
\be
R_k(z)={1\over (n\sin\al)^k}\sum_{j=0}^\infty \hat C_j(k,\al){u^j\over
  \sin^{j}\al}.\la{Rser}
\ee
The matrix elements of $C_j$ and $\hat C_j$ are bounded functions of $\al$.
Considering the remainder term in the asymptotic 
expansion of Bessel functions, we see that
$R(z)-I-\sum_{j=1}^{k-1}R_j(z)$
is given by the same series (\ref{Rser}) with different 
matrices $\hat C_j(k,\al)$ (but also bounded in $\al$). 
By analyticity, these series and their derivative w.r.t. $u$ 
converge in $U_\de$. This is also true for $\al>\al_0$.
Thus,
\be\eqalign{
R(z)=I+\sum_{j=1}^{k-1}R_j(z)+O((n\sin(\al/2))^{-k}),\\
{d\over dz}R(z)=\sum_{j=1}^{k-1}{d\over dz}R_j(z)+
O(n^{-k}\sin(\al/2)^{-k-1}), \qquad z\in U_\de,}\la{RRp}
\ee
where the remainder terms are uniform 
for $z\in U_\de$, $2s/n<\al\le\pi-\ep$, $n>s$, $s>s_0$.

Using 
(\ref{aspol},\ref{P},\ref{E},\ref{Q3}), 
we obtain for $z$ in the intersection of $ U_\de$ and
region 3 an asymptotic equivalence: 
\be\eqalign{
\chi_n^{-1}\phi_n(z)\sim\ga^n z^{n/2}\sqrt{{\pi n\over
    2}\sqrt{\om(z)}}
\left[a^{-1}I_0+aI_0'+
\sum_{j=1}^\infty\left\{
(R_{j,\,11}(z)-iR_{j,\,12}(z))a^{-1}I_0+\right.\right. \\ \left.\left.
(R_{j,\,11}(z)+iR_{j,\,12}(z))a I_0'\right\}\right],}\la{polas}
\ee
where the Bessel functions $I_0$ and $I_0'=I_1$ are taken at 
$n\sqrt{\om(z)}$. We now estimate the ``$R$'' terms and some of their 
derivatives at the point $e^{i\al}$. 
From (\ref{plem}) and (\ref{de12}) we get:
\[\eqalign{
R_{1,\,11}(e^{i\al})-iR_{1,\,12}(e^{i\al})=
{1\over 2\pi i}\int_{\partial U_\de\cup\partial\ti U_\de}
{\De_{1,\,11}(x)-i\De_{1,\,12}(x)\over x-e^{i\al}}dx=\\
{1\over 2\pi i}\oint {3a^2\over 8\sqrt{\ze}}{du\over u}+
{1\over 2i\sin\al}\mathrm{Res}_{x=e^{-i\al}}{a^2\over 8\sqrt{\ze}},}
\]
where $u=z-e^{i\al}=\epsilon e^{it}$ describes a circle of a small
radius $\epsilon$ in the positive direction. 
The expansion of $\ze$ near $e^{i\al}$ is given by
(\ref{g_as1}), and that near $e^{-i\al}$ is obtained from it by
changing the sign of $\alpha$. Expanding also $a^2(z)$ near $e^{\pm
  i\al}$ and calculating residues, we obtain 
\be
R_{1,\,11}(e^{i\al})-iR_{1,\,12}(e^{i\al})=
{3 e^{i\al/2}+e^{-i\al/2} \over 16i n\sin(\al/2)}.
\ee
Similarly, we calculate
\be
R_{1,\,11}(e^{i\al})+iR_{1,\,12}(e^{i\al})=
e^{-i\al/2}{1+e^{-i\al}-2e^{i\al}\over 3\cdot 16i n\sin(\al/2)}.
\ee

Now differentiating (\ref{plem}) w.r.t. $z$, we get
\[
R'_1(z)={1\over 2\pi i}\int_{\partial U_\de\cup\partial\ti U_\de}
{\De_1(x)\over (x-z)^2}dx.
\]
From here we obtain as above
\be
R'_{1,\,11}(e^{i\al})-iR'_{1,\,12}(e^{i\al})=
{e^{-i\al/2}\cos(\al/2)\over 16 n \sin^2(\al/2)}.
\ee

To estimate $R_2(z)$ using (\ref{plem}), we need to know 
$L(z)=R_1(z)\De_1(z)+\De_2(z)$ in the neighborhoods of $e^{\pm
  i\al}$. Here $\De_1(z)$, $\De_2(z)$ are given by (\ref{de12}) and 
for $R_1(z)$ we use the second formula in (\ref{R1}).
Then we obtain in the same way as above:
\be\eqalign{
R_{2,\,11}(e^{i\al})-iR_{2,\,12}(e^{i\al})=
{-6+7e^{i\al}-17e^{-i\al}\over 2^9 n^2\sin^2(\al/2)};\\
R_{2,\,11}(e^{i\al})+iR_{2,\,12}(e^{i\al})=
{16-9e^{i\al}+43e^{-i\al}-2e^{-2i\al}\over 3\cdot 2^9
  n^2\sin^2(\al/2)}.}
\ee

Substituting these expressions into 
(\ref{polas}), expanding $a(z)$, Bessel functions, and $\om(z)$ at
$e^{i\al}$ (see (\ref{g_as1})), and recalling (\ref{RRp}),
we obtain (\ref{phin}).
A simple calculation of (\ref{aschi}) yields (\ref{chin}). 
Taking the derivative of (\ref{polas}) at $z=e^{i\al}$ and using
the expressions for $R(e^{i\al})$ and $R'(e^{i\al})$ we complete the
proof of Theorem 2. $\Box$

\section{Proof of Theorem 1}

Consider the following weight function on the unit circle:
\[
f_\al(\th)=\cases{1, & $\al\le\th\le2\pi-\al$,\cr 0, & otherwise.},
\qquad 0<\al<\pi
\]
and the corresponding orthonormal polynomials 
$\phi_k(z,\al)=\chi_k z^k+b_{k-1}z^{k-1}+\dots+b_0$
satisfying (\ref{phi1}).
Since $f_\al(\th)=f_\al(2\pi-\th)$, the coefficients of the
polynomials $\phi_k(z,\al)$ are real. 

Associated with $f_\al(\th)$ is 
an $(n+1)\times (n+1)$ Toeplitz matrix $T_n(\al)$
whose matrix elements are as follows:
\[
(T_n(\al))_{jk}=
{1\over 2\pi}\int_0^{2\pi}e^{-i(j-k)\theta}f_{\al}(\theta)d\theta=
\cases{1-\al/\pi, & $j=k$,\cr -{\sin(\al(j-k))\over\pi(j-k)},& $j\neq
  k$},
\qquad j,k=0,1,\dots,n.
\]
Putting $\al=2s/n$ and taking the limit $n\to\infty$, we easily obtain
\be
\De(s)=\lim_{n\to\infty}\det T_{n-1}\left({2s\over n}\right).\la{tf}
\ee
If $\al$ is fixed, the large $n$ asymptotics of $\det T_n(\al)$ were
obtained by Widom \ci{Warc} (see also \ci{asto} for an alternative
derivation), namely,
\be
\det T_{n-1}(\al)=
 \cos^{n^2}(\al/2)\left(n\sin{\al\over
  2}\right)^{-1/4}2^{1/12}e^{3\zeta'(-1)}(1+o(1)),\qquad 0<\al<\pi.
\la{arc}
\ee
An idea of Dyson \ci{Dyson} was to put $\al=2s/n$ in these
asymptotics. Then, in view of (\ref{tf}), one formally obtains
the first 3 terms of (\ref{exp}). 
However, as noted in \ci{Dyson}, since the remainder 
term in (\ref{arc}) is unknown, we cannot say if this expansion 
is uniform in $\al$ and therefore cannot justify such a limit.
 
Recently, Deift \ci{D} found a formula which connects the determinants
of two different Toeplitz matrices. A proof of it we need below was
noticed by Simon. We shall
use this formula together with equation (\ref{arc}) for  
$\det T_{n-1}(\al)$ ($\al$ fixed) to obtain an expression 
for $\det T_{n-1}(2s/n)$.
The variant of Deift's formula we shall use is the following:
\be
{d\over d\al}\ln\det T_{n-1}(\al)=
-\sum_{k=0}^{n-1}{1\over 2\pi}\int_\al^{2\pi-\al}
{d\over d\al}|\phi_k(e^{i\th},\al)|^2 d\th,\qquad n=1,2,\dots\la{2det}
\ee
Indeed, it is well-known (e.g., \ci{Sz}) that the Toeplitz determinant has the
following representation in terms of the leading coefficients of
$\phi_n(z,\al)$:
\be
\det T_{n-1}(\al)=\prod_{k=0}^{n-1}\chi_k^{-2},\qquad n=1,2,\dots
\ee
Therefore 
\[
{d\over d\al}\ln\det T_{n-1}(\al)=-2\sum_{k=0}^{n-1}{\chi_k'(\al)\over
\chi_k(\al)}.
\]
On the other hand, the orthogonality of our polynomials implies
\[\eqalign{
-\sum_{k=0}^{n-1}{1\over 2\pi}\int_\al^{2\pi-\al}
{d\over d\al}|\phi_k(e^{i\th},\al)|^2 d\th=\\
-\sum_{k=0}^{n-1}{1\over 2\pi}\int_\al^{2\pi-\al}
\left(\phi_k(e^{i\th},\al)(\chi_k'(\al)e^{-ik\th}+b_{k-1}'(\al)e^{-i(k-1)\th}
+\dots)+\mathrm{\, c.c. \,}\right) d\th=
-2\sum_{k=0}^{n-1}{\chi_k'(\al)\over
\chi_k(\al)},}
\]
and (\ref{2det}) is obtained.
A corollary of it is the following 

\bl
Let $\phi_n(z,\al)$ and $T_{n-1}(\al)$ be defined as at the beginning
of the section, and
$\phi_n'(z,\al)=(d/dz)\phi_n(z,\al)$. Then for any $n=1,2,\dots$ 
\be
{d\over d\al}\ln\det T_{n-1}(\al)=
{n\over\pi}|\phi_n(e^{i\al},\al)|^2-{1\over\pi}\left\{
\phi_n(e^{-i\al},\al)e^{i\al}\phi_n'(e^{i\al},\al)+ \mathrm{
  c.c.}\right\}.\la{2det2}
\ee
\el

\noindent 
{\it Proof } From the identity
\[
{d\over d\al}\left(
\sum_{k=0}^{n-1}{1\over 2\pi}\int_\al^{2\pi-\al}
|\phi_k(e^{i\th},\al)|^2  d\th\right)={d\over d\al} n=0
\]
we obtain using (\ref{2det}),
\[
{d\over d\al}\ln\det T_{n-1}(\al)=- {1\over\pi}
\sum_{k=0}^{n-1}|\phi_k(e^{i\al},\al)|^2.
\]
As is known (e.g., \ci{Sz}), any system of orthonormal
polynomials on the circle
satisfies an analogue of the Christoffel formula which we write for $x$
and $y$ {\it on} the unit circle in the form:
\[
\sum_{k=0}^{n-1}\phi_k(x)\overline{\phi_k(y)}=
{\phi_n^*(x)\overline{\phi_n^*(y)}-\phi_n(x)\overline{\phi_n(y)}
\over 1- x{\bar y}}=
{\phi_n(x)\overline{(\phi_n(x)-\phi_n(y))}-
\phi_n^*(x)\overline{(\phi_n^*(x)-\phi_n^*(y))}
\over x\overline{(x-y)}},
\]
where $\phi_n^*(x)=x^n \overline{\phi_n(x)}$.
Letting $y\to x$ along the unit circle and noting that
$\overline{(d/dx)\phi_n^*(x)}=n x^{-n+1}\phi_n(x)-
  x^{-n+2}\phi_n'(x)$, we obtain
\be
\sum_{k=0}^{n-1}|\phi_k(x)|^2=
x\overline{\phi_n(x)}\phi_n'(x)+{1\over
  x}\phi_n(x)\overline{\phi_n'(x)}-n|\phi_n(x)|^2,\qquad |x|=1,
\ee
which, after substituting $x=e^{i\al}$ and recalling that 
the coefficients of our polynomials are real,
completes the proof. $\Box$ 

Note \cite{DIZ,TW2} that the logarithmic derivative (\ref{2det2})
satisfies a $\tau$-function version of Painlev\'e VI equation.

Now we want to integrate (\ref{2det2}) over $\al$ between $\al_1=2s/n$
and any fixed $\al_2<\pi$ for large $n$ (for $\al_2$ asymptotics
(\ref{arc}) are valid).
We therefore need to know large $n$ asymptotics of $\phi_n(z,\al)$ 
and its derivative at
the endpoint $z=e^{i\al}$ of the orthogonality arc. 
These are given by Theorem 2.
Substituting them into (\ref{2det2}), we see that the
term ${n\over\pi}|\phi_n(e^{i\al},\al)|^2$ cancels at once, and 
purely imaginary terms in the rest of the expression also disappear
(in particular, the terms with $r^3_-$ drop out).
As a result we have:
\be
{d\over d\al}\ln\det T_{n-1}(\al)=
-n^2{\sin(\al/2)\over 2\cos(\al/2)}-{\cos(\al/2)\over 8\sin(\al/2)}+
O\left({1\over n\sin^2(\al/2)}\right),
\ee 
where the remainder term is uniform for $2s/n\le\al\le\pi-\ep$, $s>s_0$,
$n>s$, $\ep>0$. 
Integrating this expression over $\al$ from $2s/n$ to $\al_2$ and
using Widom's asymptotics (\ref{arc}) for $\det T_{n-1}(\al_2)$, we
get
\be
\ln\det T_{n-1}(2s/n)=
n^2\ln\cos{s\over n}-{1\over 4}\ln n\sin{s\over n}+c_0+
O\left({1\over n\sin(s/n)}\right)+o(1)
\ee
as $n\to\infty$ with the first remainder term turning into $O(1/s)$ valid
for all $s>s_0$.
By (\ref{tf}), this equation yields (\ref{exp}).
$\Box$

\noindent
{\bf Remark} The Riemann-Hilbert problem methods allow us to calculate
asymptotics of orthogonal polynomials to arbitrary precision.
Because of the integral identity (Theorem 2b and equation (4) of
\ci{H}, see also \ci{Wsing}) that expresses Toeplitz determinants in
terms of orthogonal polynomials,
an arbitrary number of terms in asymptotics of
Toeplitz determinants could also be computed.
Thus one should be able to estimate the remainder term in (\ref{arc})
and therefore give another proof of Theorem 1.

\section{Acknowledgements} 
I am grateful to P. Deift for helpful discussions. I also thank
A. Its and H. Widom for useful comments and to 
S. Jitomirskaya for inviting me to the AMS conference at Snowbird,
where I heard \ci{D}. This work was partly supported by Sfb 288.

\end{document}